%% file: tensor-h-eigenpair.tex
\documentclass[onefignum,onetabnum]{siamart171218}

\input{ex_shared}

\ifpdf
\hypersetup{
  pdftitle={Dominant H-Eigenvectors of Tensor Kronecker Products Do Not Decouple},
  pdfauthor={Ayush Kulkarni, Charles Colley, David F. Gleich}
}
\fi

\myexternaldocument{ex_supplement}

\begin{document}

\maketitle

\begin{abstract}
  We illustrate a counterexample to an open question related to the dominant H-eigenvector of a Kronecker product of tensors. For matrices and Z-eigenvectors of tensors, the dominant eigenvector of a Kronecker product decouples into a product of eigenvectors of the tensors underlying the Kronecker product. This does not occur for H-eigenvectors and indeed, the largest H-eigenvalue can exceed the product of the H-eigenvalues of the component tensors. Beyond this general counterexample, we show this decoupling does hold in the case of diagonal tensors as well as nonnegative tensors. 
\end{abstract}



\section{Introduction}
In \citet{Colley-2023-tensor-kron}, we showed that the dominant tensor $Z$-eigenvector of a Kronecker product of tensors $\cmA$ and $\cmB$ can be expressed as the Kronecker product of the dominant tensor $Z$-eigenvectors of $\cmA$ and $\cmB$, respectively. Our result mirrors a theorem about matrix eigenvectors and Kronecker products. Here, we show that this statement is false in general for $H$-eigenvectors of tensors via a counterexample, which resolves a question that arose in the review of~\citet{Colley-2023-tensor-kron}. 

The tensor Kronecker product generalizes the Kronecker product of matrices as follows 
\[ \mbox{\includegraphics[width=\textwidth]{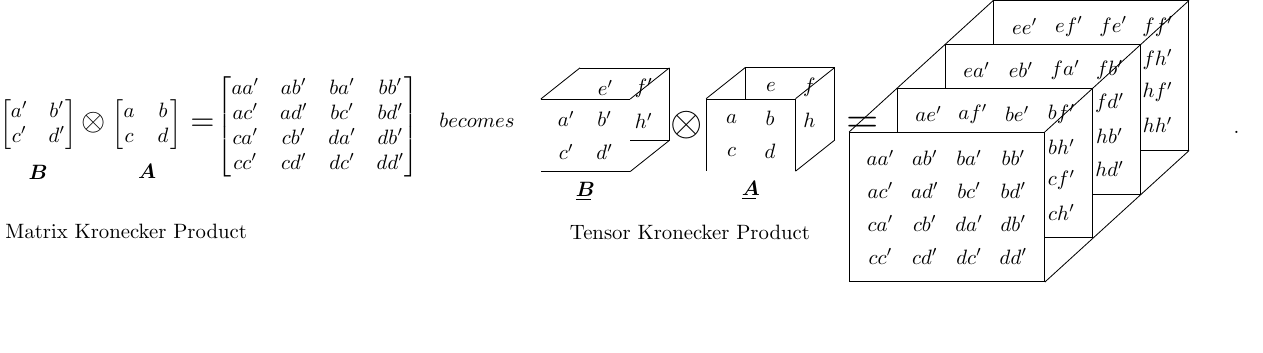}} \]
We use the increasingly standard notation where  
\[[\cmA \vx^{k-1}]_i  = \sum_{j,\ell,\ldots,s} A(i,j,\ell,\ldots, s) x_j x_\ell \cdots x_s \text{ and } [\vx^{[k-1]}]_i = x_i^{k-1}. \] 
Then a $Z$-eigenvector of a $k$-mode symmetric tensor is any solution of 
\begin{equation} \label{eq:Zeigen} \cmA \vx^{k-1} = \lambda \vx \text{ with } \normof{\vx} = 1 \text{ and  $Z$-eigenvalue }  \lambda \text{ where $\lambda \in \mathbb{C}$} \end{equation}
whereas an 
$H$-eigenvector (see note below on $\mathbb{C}$ vs $\mathbb{R}$) of a tensor is any solution of 
\begin{equation} \label{eq:Heigen} \cmA \vx^{k-1} = \lambda \vx^{[k-1]} \text{ with } \normof{\vx} = 1 \text{ and  $H$-eigenvalue } \lambda \text{ where $\lambda \in \mathbb{C}$}. \end{equation}
For more on these, see~\cite{Qi-2005-Z-eigenvalues,Lim-2005-eigenvalues}. One key difference is that for $H$-eigenvectors, the spectrum is invariant to the scale of the vector, whereas this is not the case for $Z$-eigenvectors. Note that, strictly speaking, \citet{Qi-2005-Z-eigenvalues} defines $H$-eigenvalues and eigenvectors as \emph{real} and uses the term \emph{eigenvalue and eigenvector} to represent real and complex cases. We use the terms $H$-eigenvalue or $H$-eigenvector to represent real or complex solutions of~\eqref{eq:Heigen} to clearly distinguish from the $Z$-eigenvector case. 

Let $\rho(\mA)$ be the spectral radius of a matrix, which is the largest magnitude of any eigenvalue. It is well-known that for matrices, $\rho(\mB \kron \mA) = \rho(\mA) \rho(\mB)$ and moreover, that the eigenvector that achieves the spectral radius is a Kronecker product of the eigenvectors of $\mA$ and $\mB$. 

This theorem generalizes to $Z$-eigenvectors of tensors as well. 
\begin{theorem*}[From~\citet{Colley-2023-tensor-kron}] 
		Let $\cmA$ be a symmetric, $k$-mode, $m$-dimensional tensor and $\cmB$ be a symmetric, $k$-mode, $n$-dimensional tensor. Suppose that $(\lambda_A^*, \vu^*)$ and $(\lambda_B^*, \vv^*)$ are any dominant tensor $Z$-eigenvalues and vectors of $\cmA$ and $\cmB$, respectively. Then $(\lambda_A^*\lambda_B^*, \vv^* \kron \vu^*)$ is a largest magnitude eigenpair of $\cmB \kron \cmA$.  Moreover, any Kronecker product of $Z$-eigenvectors of $\cmA$ and $\cmB$ is a $Z$-eigenvector of $\cmB \kron \cmA$. 
		\end{theorem*}

\section{Counterexamples}

We show one counterexample and one false counterexample. The first is a simple and most general case that is worked out in full detail. It has real-valued dominant eigenvalues and we are able to get analytic solutions for everything. The second false counterexample illustrates why we want to view $H$-eigenvalues as taking real or complex values. This example \emph{would} be a counterexample for a strict real-valued definition of $H$-eigenvalues but \emph{would not} be a counter-example for the more general complex-valued case we study. 

To establish our counterexamples, we need to recall a theorem from~\citet{Qi-2005-Z-eigenvalues} regarding the number of $H$-eigenvalues of a tensor. 
\begin{theorem*}[Paraphrased from \citet{Qi-2005-Z-eigenvalues}, Theorem 1, parts (a,b)] Let $\cmA$ be a symmetric tensor with $m$ modes and $n$ dimensions in each mode. (a) A number $\lambda \in \CC$ is an eigenvalue ($H$-eigenvalue over $\CC$) of $\cmA$ if and only if it is a root of the characteristic polynomial $\det(\cmA - \lambda \cmI)$ for the hyperdeterminant. (b) The number of eigenvalues ($H$-eigenvalues over $\CC$) is $n(m-1)^{n-1}$. Their product is equal to $\det(\cmA)$ (the hyperdeterminant). 
\end{theorem*}

\subsection{The simple counterexample}

The following pair of symmetric tensors is a counterexample. Let 
\begin{equation}  
\begin{array}{r@{\,}lcr@{\,}l} 
\cA(1,1,1)& = 0.3 && \cB(1,1,1) & = 0.7 \\
\cA(1,2,1)& = -0.3 &\text{ and } & \cB(1,2,1) &= -0.2 \\
\cA(1,2,2)& = 0 && \cB(1,2,2) &= -0.2 \\
\cA(2,2,2)& = 1 && \cB(2,2,2) &= -0.8 
\end{array}
\end{equation}
where the other entries are filled by symmetry.
We used a number of techniques to compute and verify the $H$-eigenvalues. For these $2 \times 2 \times 2$ tensors, we can analytically determine the $H$-eigenvalues through the polynomial from the hyperdeterminant of $(\cmA - \lambda \cmI)$. For $\cmA$, the polynomial is $-9/500 - (84/125) \lambda + (229/100) \lambda^2 - (13/5)\lambda^3 + \lambda^4$. For $\cmB$ we have the polynomial $279/625 - (9/125)\lambda - (27/20) \lambda^2 + (1/5) \lambda^3 + \lambda^4$.\footnote{To compute these expressions, we implemented the formula for the hyperdeterminant from Wikipedia and substituted in the entries of $\cmA - \lambda \cmI$ in a symbolic computing package. As further evidence that these are correct and complete, we note that both satisfy the property that $\det(\cmA) = \prod \lambda(\cmA) = -9/500$ and $\det(\cmB) = \prod \lambda(\cmB) = 279/625$ as shown in Qi's theorem. Moreover, we have $4 = 2(3-1)$ distinct values following the same theorem.} Numerically, the roots of these polynomials are
\[ \lambda(\cmA) = \{ 1, 0.812327806563 \pm 0.264915863899i, -0.024655613126 \} 
    \]
\[ \lambda(\cmB) = \{ -0.70932967445, 
        0.771909217754 \pm 0.111810698762i,
        -1.034488761057 \}. 
    \]    
The eigenvectors associated with each tensor are given by the following expression
\[ \begin{aligned}
    &  \text{for $\cmA$ } \quad \vv(\lambda) = \bmat{ 1 \\ \tfrac{1}{2} - \tfrac{5}{3} \lambda} \text{ except for $\lambda=1$ when $\vv = \bmat{0 \\ 1}$} \text{ and } \\ 
& \text{for $\cmB$ } \quad \vv(\lambda) = \bmat{ 1 \\ \frac{30 - 5 \lambda - 50 \lambda^2}{20 \lambda + 12}}. \end{aligned} \]
These expressions arise from scaling the eigenvector $\vv$ as $\vv = \sbmat{1 \\ t}$ (assuming $v_1 \not= 0$) and then solving the $H$-eigenvector equation. For the tensor $\cmA$, this simplifies due to the value of $0$ in the tensor and it becomes a simple linear equation to solve. For the tensor $\cmB$, this expansion gives two quadratic equations in $t$ we need to satisfy. We solve one for $t^2$ and substitute that value of $t^2$ into the other to get a linear function for $t$. 


We also found that \verb#HomotopyContinuation.jl#~\citep{HomotopyContinuation.jl} will solve the system of equations for the $H$-eigenvalues reliably. This software gives the same results.

Consequently, the largest magnitude $H$-eigenvalue we would  expect for the Kronecker product of $\cmB$ and $\cmA$ (if the theorem held for $H$-eigenvectors) is $1.034488761057$.

The tensor Kronecker product $\cmC = \cmB \kron \cmA$ has entries 
\[ \begin{array}{c}
C(1,1,1) = 0.21 \\
C(1,1,2) = -0.21\\
C(1,1,3) = -0.06\\
C(1,1,4) = 0.06\\
C(1,2,2) = 0.0
\end{array}
\begin{array}{c}
C(1,2,3) = 0.06\\
C(1,2,4) = -0.0\\
C(1,3,3) = -0.06\\
C(1,3,4) = 0.06\\
C(1,4,4) = -0.0\\
\end{array}
\begin{array}{c}
C(2,2,2) = 0.7\\
C(2,2,3) = -0.0\\
C(2,2,4) = -0.2\\
C(2,3,3) = 0.06\\
C(2,3,4) = -0.0
\end{array}
\begin{array}{c}
C(2,4,4) = -0.2\\
C(3,3,3) = -0.24\\
C(3,3,4) = 0.24\\
C(3,4,4) = -0.0\\
C(4,4,4) = -0.8.
\end{array}
\]
The \texttt{HomotopyContinuation.jl} software gives the following $H$-eigenvector \[ \vx = \bmat{0.099076279319 \\ 0.427548807059 \\ 
         -0.034228101784 \\ 0.89789439552} \text{ with $H$-eigenvalue $-1.035240007957$.} \] One can easily verify this expression for an H-eigenvalue and vector pair  holds computationally to machine precision. Since this value is larger (in magnitude) than $1.0349$, we have a counterexample. If we reshape this vector $\vx$ to a $2 \times 2$ matrix, then the matrix is rank $2$ (singular values 0.995, 0.105), which shows that it is \emph{not} the Kronecker product of two smaller vectors. 

\subsection{A counterexample over the real numbers that is not a counterexample over complex numbers}

After some initial experiments, we eventually settled on \verb#HomotopyContinuation.jl# to solve for both real and complex $H$-eigenvalues, as well as the hyperdeterminant formulation. Our initial studies, however, used  the AReigST software~\citep{cui2014all} for real $H$-eigenvalues. This led to the following counterexample over real-valued $H$-eigenvalues. We use $2 \times 2 \times 2 \times 2$ tensors
\[ \begin{aligned}
\underline{A}(1,1,1,1) &= -0.5 & \underline{B}(1,1,1,1) &= -0.3 \\
\underline{A}(1,1,1,2) &= 0.6  & \underline{B}(1,1,1,2) &= 1.0 \\
\underline{A}(1,1,2,2) &= 0.2  & \underline{B}(1,1,2,2) &= 0.1 \\
\underline{A}(1,2,2,2) &= 0.1  & \underline{B}(1,2,2,2) &= -1.0 \\
\underline{A}(2,2,2,2) &= -1.7 & \underline{B}(2,2,2,2) &= -0.9
\end{aligned}
\] 
where the other entries are filled by symmetry.

Like before, we compute the characteristic polynomials symbolically.  For $\cmA$, the polynomial is $-8510413/500000 - (98769/3125)\lambda - (8868/625)\lambda^2 + (1442/125)\lambda^3 + (1563/100)\lambda^4 + (33/5)\lambda^5 + \lambda^6$. For $\cmB$ we have the polynomial $12771187/250000 - (159039/3125)\lambda - (33027/625)\lambda^2 + (531/25)\lambda^3 + (339/20)\lambda^4 + (18/5)\lambda^5 + \lambda^6$. Numerically, the roots of these polynomials are:
\[ \begin{aligned}
    \lambda(\cmA) & = \{ 1.219948829489, \; -1.911703027597, \; -1.204515922037 \pm 0.965878189355i, \; \\
        & \quad ~-1.749606978909 \pm 0.022456993661i \}. \\
 \lambda(\cmB) & = \{ -2.142282147691, \; 0.725893941027, \; -1.605950142158, \; 1.225333008282, \; \\ 
 & \quad  ~-0.901497329729 \pm 3.985110194539i \}. 
\end{aligned}
    \]    
Using software that restricts search to real numbers, one would identify the spectral radii as the largest real eigenvalues: $\rho_{\mathbb{R}}(\cmA) = 1.911703027597$ and $\rho_{\mathbb{R}}(\cmB) = 2.142282147691$. If the theorem $\rho(\cmB \otimes \cmA) = \rho(\cmB)\rho(\cmA)$ applied strictly to real spectra, we would expect the largest real eigenvalue of the Kronecker product $\cmC = \cmB \otimes \cmA$ to be
\[ 1.911703027597 \times 2.142282147691 = 4.09540726771. \]
The \texttt{HomotopyContinuation.jl} software finds 64 real solutions for the system associated with $\cmC$. The largest magnitude among these real $H$-eigenvalues is approximately $4.149659269146$, corresponding to the real eigenvector
\[ \vx = \begin{bmatrix} 0.45715066048 \\ 0.594404973209 \\ 0.580380963943 \\ -0.317575090551 \end{bmatrix}. \]
Since $4.149659269146 > 4.09540726771$, this appears to be a counterexample. However, in this example, the true spectral radius of $\cmB$ is determined by its complex eigenvalues. The complex pair $-0.901497329729 \pm 3.985110194539i$ in $\lambda(\cmB)$ has a magnitude of approximately $4.085804779738$. When we calculate the product using the true spectral radii, we obtain
\[ \rho(\cmB) \rho(\cmA) = 1.911703027597 \times 4.085804779738 = 7.810845367597. \]
Solving the system for $\cmC$ over the complex numbers reveals solutions with exactly this magnitude. Specifically, \texttt{HomotopyContinuation.jl} identifies the complex $H$-eigenvalue:
\[ \lambda(\cmC) = 1.723395174615 + 7.61834722421i, \quad \text{with } |\lambda(\cmC)| = 7.810845367597. \]
Thus, the equality $\rho(\cmC) = \rho(\cmB)\rho(\cmA)$ holds when the full complex spectrum is considered.



\section{Discussion and special cases}
We initially thought this theorem would also be true for $H$-eigenvalues because Qi found their properties more similar to traditional eigenvalues. Consequently, it was a surprise to find these counterexamples. There are special cases where it does hold, however. We begin with the simplest case. 

\begin{theorem}
    Let $\cmA$ and $\cmB$ be diagonal tensors. Then $\cmB \kron \cmA$ is also diagonal. Moreover, the $H$-eigenvalues of $\cmB \kron \cmA$ are simply the products of the $H$-eigenvalues of $\cmB$ and $\cmA$. 
\end{theorem}
This theorem holds because the $H$-eigenvalues of a diagonal tensor are simply the diagonal elements, which is shown in part (c) of Theorem~1 in \citet{Qi-2005-Z-eigenvalues}.

A more interesting case is nonnegative tensors. We recall a Perron-Frobenius type result about H-eigenvalues and eigenvectors for tensors. 
\begin{theorem*}[Theorem 1.4 from \citet{ChangPearsonZhang}]
If $\cmA$ is an irreducible, nonnegative tensor of order $m$ and dimension $n$, then there exists $\lambda_0 \ge 0$ and a nonnegative vector $\vx_0 \not= 0$ such that (i) $\lambda_0$ is an $H$-eigenvalue; (ii) $\vx_0 > 0$, i.e. all components of $\vx_0$ are positive; (iii) if $\lambda$ is an eigenvalue with a nonnegative eigenvector, then $\lambda=\lambda_0$; (iv) if $\lambda$ is any eigenvalue of $\cmA$, then $|\lambda|\le \lambda_0$.
\end{theorem*}
\begin{theorem} \label{thm:nonneg}
    Let $\cmA$ and $\cmB$ be nonnegative, irreducible tensors where the Kronecker product $\cmB \kron \cmA$ is also irreducible. Then let $(\lambda_A, \vu^*)$ and $(\lambda_B, \vv^*)$ be the dominant $H$-eigenvalues and vectors of $\cmA$ and $\cmB$, respectively. Then $(\lambda_A \lambda_B, \vv^* \kron \vu^*)$ is the largest magnitude H-eigenpair of $\cmB \kron \cmA$. 
\end{theorem}

\begin{proof}
    
    Since $\cmA,\cmB$ are nonnegative and irreducible, \citet[Theorem 1.4, above]{ChangPearsonZhang} guarantees dominant $H$-eigenpairs
  $(\lambda_A,\vu^*)$, $(\lambda_B,\vv^*)$ with $\lambda_A=\rho(\cmA)>0$, $\lambda_B=\rho(\cmB)>0$
  and $\vu^*>0$, $\vv^*>0$.  

  Let $\cmC=\cmB\kron\cmA$. Since $\cmA$ and $\cmB$ are both nonnegative, then so is $\cmC$.   We know that $\vv^* \kron \vu^*$ is an $H$-eigenvector of $\cmC$ with $H$-eigenvalue $\lambda_A \lambda_B$~\cite[by $H$-eigenvector analysis of TKP Property 3, TKP Property 7]{Ragnarsson-Torbergsen2012Structured} (or also \citealt[Property 2.23]{Pickard2024}). Since $\vz = \vv^* \kron \vu^* > 0$ and $\cmC$ is irreducible  by assumption, using \citet[Theorem 1.4, part iii and iv, above]{ChangPearsonZhang} again shows that $\vz$ is an eigenvector with $H$-eigenvalue equal to the spectral radius. 
\end{proof}

Note that it is possible that $\cmA$ and $\cmB$ are irreducible, but $\cmC$ is not, so we require that assumption in the theorem. A sufficient condition to ensure this is $\cmA, \cmB > 0$ in which case $\cmC >0$ as well. 




 
 Finding other special cases remains, we believe, an interesting study.


\section*{Acknowledgments}
We used modern LLMs extensively to develop tutorials to explain the ideas among the author team, write codes, test solutions, and develop LaTeX expressions from images and sketches. The polynomial expressions for the eigenvalues were originally suggested by the LLM tools and then recomputed and verified by us. It was really quite remarkable to enter in a set of floating point eigenvalues and have the LLM output a simple polynomial that reproduces the roots exactly. These ideas led to the exact formulas using the hyperdeterminant. The LLM tools also suggested a simplified proof of Theorem~\ref{thm:nonneg} when given the statement of Theorem 1.4 from Chang, Pearson, and Zhang; our initial thoughts used an algorithmic idea based on the behavior of the power method to establish nonnegativity. We are also grateful to an anonymous referee who suggested we revisit a nonnegative ``counterexample'' mentioned in a previous revision. This ``counterexample'' turned out to be incorrect due to repeated solutions from the \texttt{HomotopyContinuation.jl} code. 

\bibliographystyle{dgleich-bib3}
\bibliography{references}

\end{document}

%% file: ex_shared.tex
\usepackage{lipsum}
\usepackage{amsfonts}
\usepackage{graphicx}
\usepackage{epstopdf}
\usepackage{algorithmic}
\ifpdf
  \DeclareGraphicsExtensions{.eps,.pdf,.png,.jpg}
\else
  \DeclareGraphicsExtensions{.eps}
\fi

\usepackage{natbib}
\usepackage{dgleich-math}
\usepackage{dgleich-tensor}

\newsiamremark{remark}{Remark}
\newsiamremark{hypothesis}{Hypothesis}
\crefname{hypothesis}{Hypothesis}{Hypotheses}
\newsiamthm{claim}{Claim}

\headers{H-Eigenpairs and Tensor Kronecker Products}{A.~Kulkarni, C.~Colley, D.~F.Gleich}

\title{Dominant H-Eigenpairs of Tensor Kronecker Products do not decouple\thanks{This work was funded by the DOE DE-SC0023162 Sparstitute MMICC center.}}

\author{Ayush Kulkarni\thanks{Metea Valley High School, Aurora, IL 
  (\email{ayushkulkarni714@gmail.com}).}
\and Charles Colley\thanks{Department of Computer Science, Purdue University, West Lafayette, IN 
  (\email{ccolley@purdue.edu}, \email{dgleich@purdue.edu}).}
\and David F.~Gleich\footnotemark[3]}

\usepackage{amsopn}

\makeatletter
\newcommand*{\addFileDependency}[1]{
  \typeout{(#1)}
  \@addtofilelist{#1}
  \IfFileExists{#1}{}{\typeout{No file #1.}}
}
\makeatother

\newcommand*{\myexternaldocument}[1]{%
    \externaldocument{#1}%
    \addFileDependency{#1.tex}%
    \addFileDependency{#1.aux}%
}